\documentclass{amsart}
\usepackage[latin1]{inputenc}

\usepackage{amssymb,amsmath}
\usepackage{amsthm}
\usepackage{paralist}
\usepackage{url}

\newtheorem{theorem}{Theorem}
\newtheorem{proposition}[theorem]{Proposition}
\newtheorem{lemma}[theorem]{Lemma}
\newtheorem{corollary}[theorem]{Corollary}

\newtheorem{problem}{Problem}
\newtheorem{fact}[theorem]{Fact}
\theoremstyle{remark}
\newtheorem{example}{Example}
\newtheorem{claim}{Claim}

\newcommand{\card}[1]{\ensuremath{\lvert{#1}\rvert}} 
\DeclareMathOperator{\ess}{ess} 
\DeclareMathOperator{\Ess}{Ess} 
\DeclareMathOperator{\gap}{gap} 

\begin{document}

\title[Join-irreducible Boolean functions]{Join-irreducible Boolean functions}

\footnotetext[1]{Institut Sup\'erieur des Technologies M\'edicales de Tunis, $^2$University of Luxembourg, $^3$Universit\'e Claude-Bernard Lyon1 and The University of Calgary, Calgary, Alberta, Canada}

\author[Moncef Bouaziz]{Moncef Bouaziz$^1$}
\address{Institut Sup\'erieur des Technologies M\'edicales de Tunis,
 9 avenue Dr Zouhair Essafi,
 1006 Tunis, Tunisie}
\email{Moncef.Bouaziz@istmt.rnu.tn}
\author[Miguel Couceiro]{Miguel Couceiro$^2$}
\address{Department of Mathematics, 
University of Luxembourg, 
162a, avenue de la Fa\"iencerie,
L-1511 Luxembourg}
\email{miguel.couceiro@uni.lu}
\author[Maurice Pouzet]{Maurice Pouzet$^3$}\thanks{The research of the first and last author has been supported by    CMCU Franco-Tunisien "Outils math\'ematiques pour l'informatique". }
\address {ICJ, Department of Mathematics, Universit\'e Claude-Bernard Lyon1, 43 Bd 11 Novembre 1918, 68622 Villeurbanne Cedex,
 France, Department of Mathematics and Statistics, The University of Calgary, 2500 University Drive NW, Calgary, Alberta, Canada T2N 1N4}
\email{pouzet@univ-lyon1.fr, mpouzet@ucalgary.ca}

\keywords{Boolean function, minor quasi-order, hypergraph, designs, Steiner systems, monomorphy.}
\subjclass[2000]{Combinatorics (05C75), (05C65), (05B05), (05B07), Order, lattices, ordered algebraic structures (06A07), (06E30), Information and communications, circuits (94C10)}

\date{\today}

\begin{abstract}  This paper is a contribution to the study of a quasi-order on the set $\Omega$ of Boolean functions, the \emph{simple minor} quasi-order. We look at the join-irreducible members of the resulting poset $\tilde{\Omega}$. Using a two-way correspondence between Boolean functions and hypergraphs, join-irreducibility translates into a combinatorial property of hypergraphs. We observe that among Steiner systems, those which yield join-irreducible members of $\tilde{\Omega}$ are the $-2$-monomorphic Steiner systems. We also describe the graphs which correspond to join-irreducible members of $\tilde{\Omega}$.
\end{abstract}
\maketitle

\section{Introduction}
Two approaches to Boolean function definability have been considered in recent years; one in terms of
 functional equations \cite{EFHH}, and one other in terms of relational constraints \cite{Pippenger}. 
 As it turned out, these two approaches have the same expressive power in the sense that they specify exactly the same classes (or properties) of Boolean functions. The characterization of these classes was first obtained by Ekin, Foldes, Hammer and Hellerstein \cite{EFHH} who showed that equational classes of  Boolean functions can be completely described in terms of a quasi-ordering $\leq$ of the set $\Omega$ of all Boolean functions, called \emph{identification minor} in \cite{EFHH, H}, \emph{simple minor} in \cite {Pippenger, CP, CL1, CL2}, \emph{subfunction} in \cite {Z}, and \emph{simple variable substitution} in \cite{CF1}.
 This quasi-order can be described as follows: for $f,g\in \Omega$, $g\leq  f $ if $g$ can be obtained from $f$ by identification of variables, permutation of variables,  and  addition or deletion of dummy variables. As shown in \cite{EFHH}, equational classes of Boolean functions
 coincide exactly with the initial segments ${\downarrow }K=\{g\in \Omega: g\leq f, \textrm{ for some }f\in K\}$ of this quasi-order, or equivalently, they  correspond to antichains $A$ of Boolean functions in the sense that they constitute sets of the form $\Omega\setminus {\uparrow }A$. Moreover,  those equational classes definable by finitely many equations  correspond to finite antichains of Boolean functions. Since then, several investigations have appeared in this direction, to mention a few, see \cite{CF1, CF2, FP, Pippenger, Po}. 
 
 The importance of Boolean function definability led to a greater emphasis on this quasi-ordering $\leq$ \cite{CP,CL1,CL2,CL3}.
As any quasi-order, the simple minor relation $\leq$ induces a partial order $\sqsubseteq$ on the set $\tilde{\Omega}$ made of equivalence classes of Boolean functions. Several properties of the resulting poset $(\tilde{\Omega},\sqsubseteq)$ were established in \cite{CP}. In particular,  it was shown that this poset is as complex as $([\omega]^{<\omega}, \subseteq)$, the poset made of finite subsets of integers and ordered by inclusion, in the sense that each is  embeddable in the other. 


 In this paper we are interested in determining the join-irreducible members of the poset $(\tilde{\Omega},\sqsubseteq)$, that is, those equivalence classes having a unique lower cover in $(\tilde{\Omega},\sqsubseteq)$. Rather than taking a direct approach, we attack this problem by looking at hypergraphs. As it is well-known, every Boolean function can be represented by a unique multilinear polynomial over the two-element field $ GF(2)$, that is, a polynomial in which each variable has degree at most one (see Zhegalkin \cite{Zhegalkin}). This  polynomial representation  of Boolean functions allows the two-way correspondence between Boolean functions and hypergraphs. 

 For any hypergraph $\mathcal{H} = (V, \mathcal{E})$ we associate the multilinear polynomial $P_{\mathcal{H}}\in GF(2)[X]$, where $X=(x_i: i\in V)$,  given by $P_{\mathcal{H}}=\underset{E\in \mathcal{E}}{\sum} \underset{i\in E}{\prod}x_i$. In fact, every multilinear polynomial $P\in GF(2)[X]$ is of the form $P=P_{\mathcal{H}}$ where $\mathcal{H} = (V, \mathcal{E})$ and $\mathcal{E}$ is the set of hyperedges corresponding to the monomials of $P$. The  simple minor relation translates into the realm of hypergraphs through the notion of quotient map. Say that a map $h'\colon V'\to V$ is a \emph{quotient map} from $\mathcal{H}' = (V', \mathcal{E}')$ to $\mathcal{H} = (V, \mathcal{E})$ if for every $E\subseteq V$, $E\in \mathcal{E}$ if and only if $\card{\{E'\in \mathcal{E}':{h'}(E')=E}$ is odd.
For two hypergraphs $\mathcal{H}'$ and $\mathcal{H}$, set $\mathcal{H}\preceq \mathcal{H}'$ if there is a quotient map from $\mathcal{H}' $ to $\mathcal{H}$. As we are going to see $\preceq$ constitutes a quasi-order between hypergraphs and two hypergraphs are related by $\preceq$ if and only if the corresponding Boolean functions are related by $\leq$ (see Lemma \ref{quasi} and Theorem \ref{correspondence}, resp.).  

The fact that a Boolean function corresponds to a join-irreducible of the poset $(\tilde{\Omega},\sqsubseteq)$ translates into a combinatorial property of the corresponding hypergraph. A description of  all hypergraphs satisfying this property eludes us.  
However, among these hypergraphs some have been intensively studied for other purposes. The basic examples are the non-trivial hypergraphs whose automorphism group is $2$-set transitive. We show that Steiner systems which yield join-irreducible members of the poset $(\tilde{\Omega},\sqsubseteq)$ are exactly those which are $-2$-\emph{monomorphic} in the sense that the induced hypergraphs obtained by deleting any pair of two distinct vertices are isomorphic (Theorem \ref {thm:steiner}). Among Steiner triple systems those with a flag-transitive automorphism group enjoy this property. We do not know if there are other.
We also describe those graphs corresponding to join-irreducible members of $(\tilde{\Omega},\sqsubseteq)$ (Theorem \ref{thm:jigraphs}). In doing so, we show that all the lower covers of each Boolean  function $f$ have the same essential arity (Theorem \ref{theorem:lowercover0}). By a result of A.Salomaa (\cite{Salomaa}, Theorem 4, p.7) it follows that this essential arity is either $\ess f-1$ or $\ess f-2$, where $\ess f$ is the essential arity of $f$. In the latter case, $f$ has (up to equivalence) a unique lower cover. This follows from Theorem \ref{booleangap} (first shown in \cite{CL1}) which provides an explicit description of those functions whose arity gap is two. 

Some of the results in this paper were presented at the ROGICS'08 conference May 12-17,  Mahdia (Tunisia) \cite{BCP}. The authors would like to express their gratitude to the organizers, Professors Y.Boudabbous and  N.Zaguia. 

 \section{Boolean functions}
A \emph{Boolean function} is simply a mapping $f\colon \{0,1\}^n\to \{0,1\}$, where $n\geq 1$. The integer $n$ is called the \emph{arity} of $f$. As simple examples of  Boolean functions we have the  \emph{projections}, i.e., mappings $(a_1,\ldots ,a_n)\mapsto a_i$, for $1\leq i\leq n$ and $a_1,\ldots ,a_m\in \{0,1\}$, and which we also refer to as \emph{variables}. 
For each $n\geq 1$, we denote by $\Omega^{(n)}=\{0,1\}^{\{0,1\}^n}$ the set of all $n$-ary Boolean functions and we denote by $\Omega=\bigcup_{n\geq 1}\Omega^{(n)}$ the set of all Boolean functions.

Let $ GF(2)$ be the two-element field $\{0,1\}$ and let $GF(2)[x_1,\ldots ,x_n]$ be the commutative ring of polynomials in the indeterminates $x_1,\ldots ,x_n$. 
To each polynomial $P\in GF(2)[x_1,\ldots ,x_n]$ corresponds an $n$-ary Boolean function $f_P \colon \{0,1\}^n\to \{0,1\}$ which is given as the evaluation of $P$, that is, for every
$(a_1,\ldots ,a_n)\in \{0,1\}^n$, $f_P(a_1,\ldots ,a_n)=P(a_1,\ldots ,a_n)$. The function $f_P$ is said to be \emph{represented} by $P$.
 As it is well-known every Boolean function can be represented in this way. In fact:

\begin{theorem}\label{Zhegalkin1} Every Boolean function $f\colon \{0,1\}^n\to \{0,1\}$, $n\geq 1$, is uniquely represented by a multilinear polynomial $P\in GF(2)[x_1,\ldots ,x_n]$ in which  each variable has degree at most one.
\end{theorem}
The multilinear polynomial $P$ is called \emph{Zhegalkin} (or \emph{Reed--Muller}) polynomial of $f$ \cite{Muller,Reed,Zhegalkin}.

A variable $x_i$ is  an \emph{essential} variable of a Boolean function $f$ if $f$ \emph{depends} on its $i$-th argument, that is, if there are $a_1,\ldots ,a_{i-1},a_{i+1}, \, \ldots ,a_{n}\in\{0,1\}$ such that
 the unary function $x \mapsto f(a_1,\ldots ,a_{i-1},x,a_{i+1},\ldots ,a_{n})$ is nonconstant. 
By \emph{essential arity} of a function $f\in\Omega^{(n)}$, denoted $\ess f$, we simply mean the number of its essential variables. 
 For instance, constant functions are exactly those functions with essential arity $0$. Functions with essential arity $1$ are either projections or negated projections.
From Theorem~\ref{Zhegalkin1} we have the following corollary.
\begin{corollary}\label{Zhegalkin2}
 A variable $x_i$ is essential in $f\in \Omega^{(n)}$ if and only if $x_i$ appears in the Zhegalkin polynomial of $f$. In particular, $\ess f$ is the number of variables appearing in the Zhegalkin polynomial of $f$.
\end{corollary}

\subsection{Simple minors of Boolean functions}

A Boolean function $g \in\Omega^{(m)}$ is said to be a \emph{simple minor} of a Boolean function $f\in \Omega^{(n)}$
if there is a mapping $\sigma \colon \{1, \ldots, n\} \to \{1, \ldots, m\}$ such that
$$g(a_1, \ldots, a_m) = f(a_{\sigma(1)}, \ldots, a_{\sigma(n)}),$$
for every $a_1, \ldots, a_m\in \{0,1\}$.
If $\sigma$ is not injective, then we speak of \emph{identification of variables.} If $\sigma$ is not surjective, then we speak of \emph{addition of inessential variables.} If $\sigma$ is a bijection, then we speak of \emph{permutation of variables.} As it is easy to verify, these Mal'cev operations are sufficient to completely describe the simple minor relation.

\begin{fact} 
 The simple minor relation between Boolean functions is a quasi-order. 
\end{fact}

Let $\leq$ denote the simple minor relation on the set $\Omega $ of all Boolean functions. If $g \leq f$ and $f \leq g$, then we say that $f$ and $g$ are \emph{equivalent,} denoted $f \equiv g$. The equivalence class of $f$ is denoted by $\tilde{f}$. If $g \leq f$ but $f \not\leq g$, then we use the notation $g < f$. The \emph{arity gap} of $f$, denoted $\gap f$, is defined by $\gap f=min\{\ess f - \ess g: g<f\}$. Note that equivalent functions may differ in arity, but not in essential arity nor in arity gap.

\begin{fact}\label{essfact} If $g \leq f$, then $\ess g \leq \ess f$, with equality if and only if $g \equiv f$.
\end{fact}

By Corollary~\ref{Zhegalkin2}, in the case of polynomial expressions, to describe the simple minor relation we only need to consider identification and permutation of essential variables, since the operation of addition of inessential variables produces the same polynomial representations. Moreover, from Fact \ref{essfact} it follows that the strict minors of a given function $f$ have Zhegalkin polynomials with strictly less variables, and that the Zhegalkin polynomials of functions equivalent to $f$ are obtained from the Zhegalkin polynomial of $f$ by permutation of its variables. 
For further developements see \cite{CL3}.

Let $(\tilde{\Omega}, \sqsubseteq)$ denote the poset made of equivalence classes of Boolean functions associated with the simple minor relation, that is, $\tilde{\Omega}=\Omega /\equiv$ together with the partial order $\sqsubseteq$ given by $\tilde{g}\sqsubseteq \tilde{f}$ if and only if $g\leq f$. Several properties of this poset were established in \cite{CP}. For example, Fact \ref{essfact} implies that each principal initial segment ${\downarrow \tilde{f}}=\{\tilde{g}:\tilde{g}\sqsubseteq \tilde{f}\}$ is finite. This means that $(\tilde{\Omega}, \sqsubseteq)$ decomposes into levels $\tilde{\Omega}_{0}, \dots  ,\tilde{\Omega}_{n}, \dots$, where $\tilde{\Omega}_{n}$ is the set of minimal elements of $\tilde{\Omega}\setminus \bigcup \{ \tilde{\Omega}_{m}: m<n\}$. For instance, the first level $\tilde{\Omega}_{0}$ comprises four equivalence classes, namely, those of constant $0$ and $1$ functions, and those of projections and negated projections. These four classes induce a partition of $(\tilde{\Omega}, \sqsubseteq)$ into four different blocks with no comparabilities in between them. These facts were observed in \cite{CP} where it was shown that each level of $(\tilde{\Omega}, \sqsubseteq)$ is finite (\cite{CP}, Corollary 1, p.75).
The latter is entailed by the following result by A.Salomaa \cite{Salomaa}.

\begin{theorem}\label{Salomaa0}
 The arity gap of any Boolean function is at most $2$.
\end{theorem}

The description of those Boolean functions with arity gap 2 is given below.

\begin{theorem}\label{booleangap}\emph{(In \cite{CL1}:)}
Let $f \colon \{0, 1\}^n \to \{0, 1\}$ be a Boolean function with at least two essential variables. Then the arity gap of $f$ is two if and only if it is equivalent to one of the following functions:
\begin{compactenum}
\item $x_1 + x_2 + \dots + x_m + c$ for some $m \geq 2$,
\item $x_1 x_2 + x_1 + c$,
\item $x_1 x_2 + x_1 x_3 + x_2 x_3 + c$,
\item $x_1 x_2 + x_1 x_3 + x_2 x_3 + x_1 + x_2 + c$,
\end{compactenum}
where $c \in \{0,1\}$ and where $+$ is taken modulo $2$. Otherwise the arity gap of $f$ is one.
\end{theorem}

Theorem~\ref{Zhegalkin1} allows to work with polynomials rather than Boolean functions. This approach turns out to be quite useful when studying the poset  $(\tilde{\Omega}, \sqsubseteq)$. For instance, the four equivalence classes in $\tilde{\Omega}_{0}$ mentioned above are represented by $0,1,x_1$ and $x_1+1$. As it is easy to verify, above the equivalence classes represented by the constant polynomials $0$ or $1$ we have the equivalence classes of those functions whose Zhegalkin polynomials are the sum of an even number of nonconstant monomials plus $0$ or $1$, respectively, and above the equivalence classes represented by $x_1$ or $x_1+1$ we have the equivalence classes of those functions whose Zhegalkin polynomials are the sum of an odd number of nonconstant monomials plus $0$ or $1$, respectively.

\subsection{Join-irreducible Boolean functions}

We say that an element $\tilde{f} \in \tilde{\Omega}$ is \emph{join-irreducible} if there is $\tilde{f}'\in \tilde{\Omega}$ such that 
$\tilde{f}'\sqsubset \tilde{f}$ and 
 for every $\tilde{g} \in \tilde{\Omega}$,  if $\tilde{g}\sqsubset \tilde{f}$, then $\tilde{g} \sqsubseteq \tilde{f}'$. 
Since $\tilde{\Omega}$ decomposes into levels this amounts to say that $\tilde{f}$ has a unique lower cover.
For the sake of simplicity, we say that a function $f\in \Omega$ is \emph{join-irreducible} if $\tilde{f}$ is join-irreducible. Likewise, we say that
$g$ is a \emph{lower cover} of $f$ if $\tilde{g}$ is a lower cover of $\tilde{f}$. 

To illustrate, consider the binary conjunction $x_1\wedge x_2$, the binary disjunction $x_1\vee x_2$. Both of these functions are join-irreducible  since they have, up to equivalence, a unique strict simple minor, namely, a projection. This uniqueness clearly extends to any conjunction and disjunction of $n\geq 2$ variables, showing that any of the latter functions also constitute join-irreducible functions. But this is not the case for the composite $(x_1\vee x_2) \wedge x_3\wedge x_4$. Indeed, $(x_1\vee x_2) \wedge x_3, \, x_1\wedge x_3\wedge x_4 < (x_1\vee x_2)\wedge  x_3\wedge x_4$, but $(x_1\vee x_2) \wedge x_3 \not \equiv x_1\wedge x_3\wedge x_4 $. These observations lead to the following problem.
  
\begin{problem}\label{mainproblem} Describe the join-irreducible Boolean functions.
 \end{problem}

Let $f\colon \{0,1\}^n\to \{0,1\}$ be a Boolean function and let $i, j\in \{1,\ldots,n\}$. We denote by $f_{i=j}$ the function obtained from $f$ by identifying the variable $x_{i}$ to the variable $x_{j}$ with the convention that $f_{i=i}=f$.

\begin{lemma} \label{lem:-1} If $g<f$, then there are two distinct essential variables $i$ and $j$ of $f$ such that $g\leq f_{i=j}<f$.\end{lemma}

Using Lemma \ref{lem:-1}, we see that each of the functions given in Theorem \ref{booleangap} is join-irreducible since each has, up to equivalence, a unique lower cover, namely, $c$ in case $(2)$ and $x+c$ in cases $(1),(3)$ and $(4)$. Thus to solve Problem \ref{mainproblem} we need to focus on 
those functions whose arity gap is equal to one.
Towards this problem we will make use of the following result.

\begin{theorem}\label{theorem:lowercover0} 
All lower covers of a given function $f$ have the same essential arity which is either $\ess f-1$ or $\ess f-2$. In the latter case,
$f$ is join-irreducible.
\end{theorem}

To prove Theorem~\ref{theorem:lowercover0}, we make use the following auxiliary properties of a Boolean function $f\colon \{0,1\}^n\to \{0,1\}$.

\begin{lemma} \label{lem:1} Let $i_1,i_2, t\in \{1,\ldots,n\}$ such that $i_1\not =i_2$. If $x_t$ is inessential in $f_{i_1=i_2}$, then 
$f_{t=k}\geq f_{{i_1= i_2},{t=k}}\equiv f_{i_1=i_2}$, for all $k\in \{1,\ldots,n\}$.
\end{lemma}

\begin{lemma}\label{lem:2} If $f_{i_1=k}\equiv f_{i_2=k}$ and $x_k$ is essential in $f_{i_1=k}$ and in $f_{i_2=k}$, then $x_{i_1}$ is inessential in $f_{i_2=k}$ if and only if $x_{i_2}$ is inessential in $f_{i_1=k}$.
\end{lemma}
\begin{proof}
 Suppose that $x_{i_1}$ is inessential in $f_{i_2=k}$. By Lemma~\ref{lem:1}, $f_{i_2=k,i_1=k}\equiv f_{i_2=k}$. Suppose for the sake of a contradiction that $x_{i_2}$ is essential in $f_{i_1=k}$. Then, since $x_k$ is essential in $f_{i_1=k}$ and in $f_{i_2=k}$, we have 
$f_{i_1=k,i_2=k}< f_{i_1=k}$. Hence, $f_{i_2=k}<f_{i_1=k}$ which constitutes the desired contradiction. 
\end{proof}

Note that the hypotheses of Lemma~\ref{lem:2} are satisfied by $f=x_1x_2+x_1x_3+x_2x_3$.

\begin{lemma}\label{lem:0}
Suppose that $i_1, i_2, k$ are distinct indices in $\{1,\ldots,n\}$ such that  
$x_{i_1}$ is inessential in $f_{i_2=k}$ and $x_{i_2}$ is inessential in $f_{i_1=k}$. If $x_{i_1}$ is inessential in $f_{i_2=i_1}$, then $x_{i_1}$ is inessential in $f$.
\end{lemma}

\begin{proof}
Without loss of generality, suppose that $i_1=1, i_2=2$ and $k=3$. Let $(a_4, \dots, a_n)\in \{0,1\}^{n-3}$ and let $g\colon \{0,1\}^3\to \{0,1\}$ be the function given by $g(x_1, x_2, x_3)=f( x_1, x_2, x_3, a_4, \dots, a_n)$. We show that 
\begin{equation}\label{eq:0}
g(a_1, a_2, a_3)= g(1+a_1, a_2, a_3)
\end{equation}
 for all $(a_1, a_2, a_3)\in \{0, 1, 2 \}^3$. 
 Indeed, if $a_2=a_3$, then
\begin{equation*}
g(a_1,a_2,a_3)=g(a_1,a_3,a_3)=g(1+a_1,a_3,a_3)=g(1+a_1,a_2,a_3)
\end{equation*} since $x_1$ is inessential in $f_{2=3}$.

Now suppose that $a_2=1+a_3$. If $a_1=a_3$, then 
 \begin{equation*}\label{eq:2} 
g(a_1,a_2,a_3)=g(a_3,1+a_3,a_3)=g(a_3,a_3,a_3)=g(a_1,1+a_2,a_3)
\end{equation*} since $x_2$ is inessential in $f_{1=3}$, and we have
 \begin{equation*}
g(a_1,1+a_2,a_3)=g(a_3,a_3,a_3)=g(1+a_3,1+a_3,a_3)=g(1+a_1,a_2,a_3)
\end{equation*} because $x_1$ is inessential in $f_{2=1}$. Thus (\ref{eq:0}) holds.

 If $a_1=1+a_3$, then
\begin{equation*}\label{eq:3} 
g(a_1,a_2,a_3)=g(1+a_3,1+a_3,a_3)=g(a_3,a_3,a_3)=g(1+a_1,1+a_2,a_3)
\end{equation*} since $x_1$ is inessential in $f_{2=1}$, and we have
\begin{equation*} 
g(1+a_1,1+a_2,a_3)=g(a_3,a_3,a_3)=g(a_3,1+a_3,a_3)=g(1+a_1,a_2,a_3)
\end{equation*}because $x_2$ is inessential in $f_{1=3}$. Thus (\ref{eq:0}) holds and the proof is now complete.
\end{proof}

\begin{lemma}\label{lemma:lowercover} 
Let $g$ be a lower cover of $f$. If $ess(g)<ess(f)-1$, then $g$ is, up to equivalence, the unique lower cover of $f$. 
\end{lemma}

\begin{proof}
With no loss of generality, we may assume that every variable of $f$ is essential. Let  $G:=(V, \mathcal E)$, where $V=\{1,\ldots,n\}$ and $\mathcal E:= \{\{i,j\}\in [V]^2: f_{i=j}\equiv g\}$. Since $g$ is a lower cover of $f$, $G$ is not the empty graph. Our aim is to show that $G$ is a complete graph.
Let $A:= \{i\in V: \{i,k\}\in \mathcal E \; \text{ for all }\; k\in V\setminus \{i\}\}$. Our aim reduces to prove that $A=V$. 
\begin{claim} \label{claim:0}
Let $\{i,j\}\in \mathcal E$. Then:
\begin{enumerate}[(a)]
\item If $x_t$ is inessential in $f_{i=j}$, then $t$ belongs to $A$.
\item If $j\not \in A$ and if $x_t$ is inessential in $f_{i=j}$, for some $t\in V\setminus \{i\}$, then $t\neq j$, $\{t,j\}\in \mathcal E$ and $x_i$ is inessential in $f_{t=j}$.
\item $i$ and $j$ belong to $A$.
\end{enumerate}
\end{claim}
\begin{proof}[Proof of Claim \ref{claim:0}] 
First we prove (a). Let $k\in V\setminus \{t\}$. According to Lemma \ref{lem:1}, $f_{i=j}\leq f_{t= k}$. Since $g$ is a lower cover of $f$ and all variables of $f$ are essential, $f_{t=k}\equiv g$, that is, $\{t, k\}\in \mathcal E$. Since this holds for every $k\in V\setminus \{t\}$, we have $t\in A$.

To prove (b), suppose that $j\not \in A$. Then from (a) it follows that $x_j$ is essential in $f_{i=j}$. 
 Since $g$ is a lower cover of $f$ and $x_t$ is inessential in $f_{i=j}$, Lemma \ref{lem:1} yields $f_{i=j}=f_{t=j}$, thus $\{t,j\}\in \mathcal E$. Since $j\not \in A$, it follows from (a) that $x_j$ is essential in  $f_{t=j}$. Applying Lemma \ref{lem:2} to $f$, with $k:=j$, $i_1:= i$, $i_2:=t$, yields that $x_i$ is inessential in $f_{t=j}$.

To prove (c) suppose, for the sake of a contradiction, that $j\not \in A$. Hence, $x_j$ is essential in $f_{i=j}$. Since $ess(g)<ess(f)-1$, there is some  $t\in V\setminus \{i,j\}$ such that $x_t$ is inessential in $f_{i= j}$ and thus $t\in A$.  From (b), it follows that $x_i$ is inessential in $f_{t=j}$. By Lemma \ref{lem:0}, $x_t$ is essential in $f_{i=t}$. Since $x_i$ is inessential in $f_{t=j}$, it follows from Lemma \ref{lem:2} that $x_j$ is inessential in $f_{i=t}$. This implies $j\in A$, a contradiction. 
\end{proof} 
From (c) it follows that $A=V$ and the proof of the lemma is complete.
\end{proof}

\begin{proof}[Proof of Theorem \ref{theorem:lowercover0}]
 By Theorem \ref{Salomaa0}, the essential arity of the lower covers of $f$ are either $\ess f-1$ or $\ess f-2$. By Lemma \ref{lemma:lowercover}, if there is a lower cover with essential arity equal to $\ess f-2$, then it is, up to equivalence, unique and thus $f$ is join-irreducible. 
\end{proof}

Set $\Ess f = \{i\in \{1, \ldots, n\}: x_i \; \text {is an essential variable of}\; f\}$ and let $[\Ess f]^2$ be the set of $2$-element subsets of $\Ess f$.
For  $e=\{i,j\}, e'=\{i',j'\}\in  [\Ess f]^2$, define $e\approx e'$ if $f_{i=j}\equiv f_{i'=j'}$. Obviously, $\approx$ is an equivalence relation  on  $[\Ess f]^2$ and $\ess f_{i=j}= \ess f_{i'=j'}$ for all  $e=\{i,j\}$, $e'=\{i',j'\}$ such that $e\approx e'$. According to Lemma \ref{lem:-1},  if a Boolean function $g$ is a lower cover of $f$, then there is some $\{i,j\} \in  [\Ess f]^2$ such that $g\equiv f_{i=j}$.  From this observation, we get the following fact.
\begin{fact} \label{fact-trivial}A Boolean function $f$ is join-irreducible provided that $\ess f\geq 2$ and $[\Ess f]^2$ is an equivalence class.
\end{fact}

To provide our first criterion for join-irreducibility, let $C_f = \{\{i,j\} \in  [\Ess f]^2:  \ess f_{i=j}=\ess f- \gap f \}$.
Clearly, for each pair $\{i,j\}\in C_f$, we have that $f_{i=j}$ is a lower cover of $f$, and Theorem \ref {theorem:lowercover0} asserts that the converse also holds.
 Hence,  $C_f $ is the union of equivalence classes of pairs $e=\{i,j\}\in  [\Ess f]^2$ such that $f_{i=j}$ is a lower cover of $f$.
For example, if $\gap f=2$, then $C_f $ is an equivalence class, namely, the whole set $ [\Ess f]^2$ (to see this, use the description given in Theorem \ref{booleangap}).
 These observations yield  our first criterion for join-irreducibility.

\begin{theorem}\label{theorem:key0}
A Boolean function $f$ is join-irreducible if and only if $\ess f\geq 2$ and $C_f$ is an equivalence class of $\approx$.
 Furthermore, if $\gap f=2$, then $f$ is join-irreducible and $C_f=[\Ess f]^2$.
\end{theorem}

To make this criterion applicable, we will encode Boolean functions by hypergraphs and translate Theorem~\ref{theorem:key0}, accordingly.

\section{Boolean functions and hypergraphs}

By an \emph{hypergraph} we simply mean a pair $\mathcal{H} = (V, \mathcal{E})$ where $V$ is a finite nonempty set whose elements are called \emph{vertices}, and where $\mathcal{E}$ is a collection of subsets of $V$ called \emph{hyperedges}. We write $[V]^m$ to denote the set of $m$-element subsets of $V$.

Let $\mathcal{H} = (V, \mathcal{E})$ be an hypergraph with $n$ vertices. To such an hypergraph $\mathcal{H}$ we associate a Zhegalkin polynomial $P_{\mathcal{H}}\in GF(2)[x_i: i\in V]$ which is given by $P_{\mathcal{H}}=\underset{E\in \mathcal{E}}{\sum} \underset{i\in E}{\prod}x_i$.

\begin{example}
Let $\mathcal{H}_1 = (\{1,2,3\}, \emptyset)$,   $\mathcal{H}_2 = (\{1,2,3\},\{\{1,2\}, \emptyset\})$ and $\mathcal{H}_3 = (\{1,2,3\}, \{\{1,2\}, \{1,3\}, \{2,3\}\} )$. Then $\mathcal{P}_{\mathcal{H}_1}= 0$, $\mathcal{P}_{\mathcal{H}_2}= x_1x_2+1$ and $\mathcal{P}_{\mathcal{H}_3}= x_1x_2+x_1x_3+x_2x_3$, respectively.
\end{example}

Conversely, to each Zhegalkin polynomial $P\in GF(2)[x_1,\ldots ,x_n]$ is associated an hypergraph 
$\mathcal{H}_P = (V, \mathcal{E})$ where $V=\{1,\ldots ,n\}$ and $\mathcal{E}$ is the set of hyperedges corresponding to the monomials of $P$. By Theorem~\ref{Zhegalkin1}, we have the following.

\begin{theorem} For each Boolean function $f\colon \{0,1\}^n\to \{0,1\}$, $n\geq 1$, there is a unique hypergraph $\mathcal{H} = (V, \mathcal{E})$, $V=\{1,\ldots ,n\}$, such that $f=f_{P_{\mathcal{H}}}$.
\end{theorem}

For the sake of simplicity, let $f_{\mathcal{H}}$ denote the function $f_{P_{\mathcal{H}}}$ determined by $\mathcal{H}$.

 \subsection{Simple minors of hypergraphs}

Let $\mathcal{H}' = (V', \mathcal{E}')$ be an hypergraph and let $h'\colon V'\to V$ be a map. For each $E\subseteq V$, set ${h'}^{-1}[E]=\{E'\in \mathcal{E}': h'(E')=E\}$, where $h'(E')=\{h'(i'): i'\in E'\}$.
If $\mathcal{H} = (V, \mathcal{E})$ is an hypergraph, then the map $h'$ is said to be a \emph{quotient map from $\mathcal{H}'$ to $\mathcal{H}$}, denoted $h'\colon \mathcal{H}'\to \mathcal{H}$, if for every 
$E\subseteq V$, the following condition holds: $E\in \mathcal{E}$ if and only if the cardinality $\card{{h'}^{-1}[E]}$ is odd.
We say that an hypergraph $\mathcal{H}$ is a \emph{simple minor} of an hypergraph $\mathcal{H}'$, denoted $\mathcal{H}\preceq \mathcal{H}'$, if there is a quotient map from $\mathcal{H}'$ to $\mathcal{H}$.

To illustrate, let $\mathcal{H} = (V, \mathcal{E})$ be an hypergraph with $V=\{1,\ldots ,n\}$. Let $e=\{i,j\}$, $i,j\in V$, and fix $l_e\not \in V$. Consider the hypergraph $\mathcal{H}_{e} = (V_{e}, \mathcal{E}_{e})$ given as follows: $V_{e}=(V\setminus e)\cup \{l_e\}$ and for each $E \subseteq V_{e}$, we have $E \in \mathcal{E}_{e}$ if either
\begin{enumerate}[(i)]
 \item \label{oddcondition1} $E \in \mathcal{E}$ and $e\cap E=\emptyset$, or
\item \label{oddcondition2} $l_e\in E$ and among $(E\setminus \{l_e\})\cup e$, $(E\setminus \{l_e\})\cup \{i\}$ and $(E\setminus \{l_e\})\cup \{j\}$, either one or all of these sets belong to $\mathcal{E}$.
\end{enumerate}
Note that conditions (\ref{oddcondition1}) and (\ref{oddcondition2}) guarantee that the map $h\colon V\to V_{e}$ defined by $h(i)=h(j)=l_e$ and $h(k)=k$, for each $k\neq i,j$, constitutes a quotient map from $\mathcal{H}$ to $\mathcal{H}_{e}$, thus showing that $\mathcal{H}_{e}$ is a simple minor of $\mathcal{H}$.

 \begin{lemma}\label{quasi}
  The simple minor relation between hypergraphs is a quasi-order.
 \end{lemma}

\begin{proof}
 Let $\mathcal{H} = (V, \mathcal{E})$, $\mathcal{H}' = (V', \mathcal{E}')$ and $\mathcal{H}^{''} = (V^{''}, \mathcal{E}^{''})$ be hypergraphs such that $\mathcal{H}\preceq \mathcal{H}^{'}\preceq \mathcal{H}^{''}$. Let $h'\colon \mathcal{H}'\to \mathcal{H}$ and $h^{''}\colon \mathcal{H}^{''}\to \mathcal{H}'$ be the corresponding quotient maps. We claim that $h=h'\circ h^{''}$ is a quotient map from $\mathcal{H}^{''}$ to 
$\mathcal{H}$, i.e., $\mathcal{H}\preceq \mathcal{H}^{''}$. 

Let $E\subseteq V$. The set $h^{-1}[E]=\{E^{''}\in \mathcal{E}^{''}:h^{''}(E'')=E \}$ decomposes into two sets, namely, 
$A=\bigcup \{{h^{''}}^{-1}[E']:E'\in \mathcal{E}', h'(E')=E\}$ and 
$B=\bigcup \{{h^{''}}^{-1}[E']:E'\not\in \mathcal{E}', h'(E')=E\}$.
Now $A$ is a disjoint union of sets of odd size and $B$ is a disjoint union of sets of even size and hence, $B$ has even size. Thus the parity of $\card{h^{-1}[E]}$ is the same as the parity of $\card{A}$ which, in turn, is the same as the parity of $\card{{h'}^{-1}[E]}$. Since $E\in \mathcal{E}$ if and only if $\card{{h'}^{-1}[E]}$ is odd, the proof is now complete.
\end{proof}

A simpler proof of Lemma~\ref{quasi} can be obtained using the following construction. Given an hypergraph $\mathcal{H}' = (V', \mathcal{E}')$, a set $V$ and a map $h'\colon V'\to V$, let $\mathcal{H}_{h'} = (V, \mathcal{E}_{h'})$ where 
$\mathcal{E}_{h'}=\{E\subseteq V\colon \card{{h'}^{-1}[E]}_2=1\}$ and where $\card{{h'}^{-1}[E]}_2$ denotes the cardinality of ${h'}^{-1}[E]$ modulo 2.

\begin{lemma}\label{lemma:construction} Let $\mathcal{H}' = (V', \mathcal{E}')$ and a map $h'\colon V'\to V$, with  $V=\{1,\ldots ,n\}$ and 
$V'=\{1,\ldots ,m\}$.
Then $h'$ is a quotient map from $\mathcal{H}'$ to $\mathcal{H}_{h'}$ and 
$$f_{\mathcal{H}_{h'}}(a_1,\ldots ,a_n)=f_{\mathcal{H}'}(a_{h'(1)}, \ldots, a_{h'(m)}),$$
for all $a_1,\ldots ,a_n\in\{0,1\}$.
Moreover, if $h'$ is a quotient map from $\mathcal{H}'$ to $\mathcal{H}$, then $\mathcal{H}=\mathcal{H}_{h'}$.  
\end{lemma}
\begin{proof}
By construction, $h'$ is a quotient map from $\mathcal{H}'$ to $\mathcal{H}_{h'}$. Using the fact that $a^2=a$, for every $a\in \{0,1\}$, we have
\begin{eqnarray*}
f_{\mathcal{H}_{h'}}(a_1,\ldots ,a_n) &=& \underset{E\in \mathcal{E}_{h'}}{\sum} \underset{i\in E}{\prod}a_i
= \underset{E\in \mathcal{E}_{h'}}{\sum} \card{{h'}^{-1}[E]}_2 \underset{i\in E}{\prod}a_i\\
&=& \underset{E\in \mathcal{E}_{h'}}{\sum}\, \underset{E'\in {h'}^{-1}[E]}{\sum}\, \underset{i\in h'(E')}{\prod}a_{i}
= \underset{E'\in \mathcal{E}'}{\sum} \underset{i\in h'(E')}{\prod}a_{i}^{\card{{h'}^{-1}(i)}}\\
&=& \underset{E'\in \mathcal{E}'}{\sum} \underset{i'\in E'}{\prod}a_{h'(i')}
= f_{\mathcal{H}'}(a_{h'(1)}, \ldots, a_{h'(m)}),
\end{eqnarray*}
for every $a_1,\ldots ,a_n\in \{0,1\}$. The last claim is an immediate consequence of the construction of $\mathcal{H}_{h'}$. 
\end{proof}

The following theorem establishes the connection between the simple minor relation on Boolean functions and the simple minor relation on hypergraphs.

\begin{theorem}\label{correspondence} Let $\mathcal{H} = (V, \mathcal{E})$ and $\mathcal{H}' = (V', \mathcal{E}')$ be two hypergraphs, with  $V=\{1,\ldots ,n\}$ and 
$V'=\{1,\ldots ,m\}$, respectively. Then $\mathcal{H}\preceq \mathcal{H}'$ if and only if $f_{\mathcal{H}}\leq f_{\mathcal{H}'}$.
\end{theorem}

\begin{proof} Let $h'\colon V'\to V$. If $h'$ is a quotient from $\mathcal{H}'$ to $\mathcal{H}$, then $\mathcal{H}=\mathcal{H}_{h'}$ and thus
$$
f_{\mathcal{H}}(a_1,\ldots ,a_n)= f_{\mathcal{H}_{h'}}(a_1,\ldots ,a_n)=f_{\mathcal{H}'}(a_{h'(1)}, \ldots, a_{h'(m)}),
$$
for every $a_1,\ldots ,a_n\in \{0,1\}$, by Lemma~\ref{lemma:construction}. 
Conversely, if 
$$
f_{\mathcal{H}}(a_1,\ldots ,a_n)=f_{\mathcal{H}'}(a_{h'(1)}, \ldots, a_{h'(m)}),
$$
for every $a_1,\ldots ,a_n\in \{0,1\}$, then $f_{\mathcal{H}}= f_{\mathcal{H}_{h'}}$, by Lemma~\ref{lemma:construction}.
By uniqueness of the Zhegalkin polynomial representation, $\mathcal{H}= \mathcal{H}_{h'}$. Hence, $h'$ is a quotient map from
$\mathcal{H}'$ to $\mathcal{H}$. This completes the proof of the theorem.
\end{proof}

\subsection{Conditions for join-irreducibility}\label{Sec:3.2}

Let $\mathcal{H} = (V, \mathcal{E})$ and $\mathcal{H}' = (V', \mathcal{E}')$ be two hypergraphs. A map  $\varphi \colon V\to V'$ is said to be an \emph{isomorphism} from $\mathcal{H}$ onto $\mathcal{H}'$ if $\varphi$ is bijective and for every $E\subseteq V$, $E\in \mathcal{E}$ if and only if $\varphi(E) \in \mathcal{E}'$. Two hypergraphs  $\mathcal{H} $ and $\mathcal{H}' $ are said to be \emph{isomorphic}, denoted $ \mathcal{H} \cong \mathcal{H}' $, if there is an isomorphism $\varphi $ from $\mathcal{H}$ onto $\mathcal{H}'$. 
If $\mathcal{H} =\mathcal{H}' $, then $\varphi$ is called an \emph{automorphism} of $\mathcal{H}$. The group made of automorphisms of $\mathcal{H}$ is denoted by $Aut(\mathcal{H})$.

Let $\mathcal{H} = (V, \mathcal{E})$ be an hypergraph. We set $\check V= \bigcup \mathcal{E}$. Since the essential variables of $f_{\mathcal{H}}$ are those which appear in the Zhegalkin polynomial,  $\Ess f_{\mathcal{H}}= \check V$, and hence, $\ess f= \vert \check V\vert $.
Note that for every $e,e'\in  [\check{V}]^2$, we have  $e\approx e'$ whenever $\mathcal{H}_e\cong \mathcal{H}_{e'}$. Set $D_{\mathcal H}= \{e\in  [\check{V}]^2:  \check {V_{e}}= (\check V\setminus e)\cup l_{e}\}$. Clearly, $e\in D_{\mathcal H}$ if and only if $\ess f_{\mathcal{H}_{e}}=\ess f_{\mathcal{H}}- 1$.

Given these observations, we obtain from  Fact \ref{fact-trivial} our first  criterion for join-irreducibility.
\begin{proposition}\label{cor:1} Let $\mathcal{H} = (V, \mathcal{E})$ be an hypergraph. If $\card{ \check{V}}\geq 2$ and $D_{\mathcal H}=[\check{V}]^2$ then $f_{\mathcal{H}}$ is join-irreducible.
\end{proposition}

We may translate Theorem \ref{theorem:key0} as follows.
\begin{theorem}\label{theorem:key}
 Let $\mathcal{H} = (V, \mathcal{E})$ be an hypergraph. Then
$f_{\mathcal{H}}$ is join-irreducible if and only if $\card{ \check{V}}\geq 2$ and either $D_{\mathcal H}$ is an equivalence class or $\mathcal H=\mathcal H_P$ where $P$ is one of the polynomials given in Theorem \ref{booleangap}. \end{theorem}



In the search for hypergraphs $\mathcal{H} = (V, \mathcal{E})$ determining join-irreducible Boolean functions, Theorem~\ref{theorem:key} invites us to look at differences $\ess f_{\mathcal{H}}-\ess f_{\mathcal{H}_{e}}$, especially, at the cases when $\ess f_{\mathcal{H}}-\ess f_{\mathcal{H}_{e}}>1$. For the latter to occur, there are two possibilities:
\begin{enumerate}[(i)]
 \item the vertex $l_e$ becomes isolated and this is the case if and only if, for every $F$ disjoint from $e$, the number of $e'\subseteq V$ such that $\emptyset \neq e'\subseteq e$ and $e'\cup F \in \mathcal{E}$, is even, or
\item another vertex, say $i\in V$, becomes isolated and this is the case if and only if, for every $e'\in \mathcal{E}$, if $i\in e'$ then $e\cap e'\neq \emptyset$ and there is exactly one $e''\in \mathcal{E}$ distinct from $e'$ such that $i\in e''$ and $e'\setminus e=e''\setminus e$.
\end{enumerate}  

Proposition \ref{cor:1} reveals some interesting connections with some well-known combinatorial properties of hypergraphs. 
A group $G$ acting on a set $V$ is \emph{$2$-set transitive} if for every $e,e'\in [V]^2$, there is some $g\in G$ such that $g(e)=e'$.

\begin{corollary} \label{2transitive}Let $\mathcal{H} = (V, \mathcal{E})$ be an hypergraph. If $\card{V}\geq 2$, $\bigcup \mathcal{E}= V$ and $Aut(\mathcal{H})$ is $2$-set transitive, then $f_{\mathcal{H}}$ is join-irreducible.
\end{corollary}
\begin{proof}
 Let $\varphi \in Aut(\mathcal{H})$. Take $e\in [{V}]^2$ and let $e'=\varphi(e)\in [{V}]^2$. Consider the mapping $\overline{\varphi}:V_e\to V_{e'}$ defined by $\overline{\varphi}(l_e)=l_{e'}$ and $\overline{\varphi}(i)=\varphi (i)$ for every $i\in V_e\setminus \{l_e\}$. Clearly, $\overline{\varphi}$ constitutes the desired isomorphism from $\mathcal{H}_e$ to $\mathcal{H}_{e'}$.
\end{proof}
Proposition  \ref{cor:1} and Corollary \ref{2transitive} naturally give raise to the following questions:
\begin{problem}
 For which hypergraphs $\mathcal{H}$
\begin{enumerate}[(i)]
 \item $Aut(\mathcal{H})$ is $2$-set transitive?
\item $\mathcal{H}_e\cong \mathcal{H}_{e'}$, for every $e,e'\in [{V}]^2$?
\end{enumerate}
\end{problem}
\subsection{Steiner Systems}

Let $\mathcal{H}=(V, \mathcal{E})$ be an hypergraph. We say that the hypergrapph $\mathcal{H}$ is a \emph{$2-(n,k,\lambda)$ design} if $\card{V}=n$, $\mathcal{E}\subseteq [V]^k$, and for every $e\in [V]^2$, $\card{\{E\in \mathcal{E}: e\subseteq E\}}=\lambda$.
 If $\lambda =1$, then we say that $\mathcal{H}$ is a \emph{Steiner system} and, in addition, if $k=3$, then we say that $\mathcal{H}$ is a \emph{Steiner triple system}.

For each $e\in [V]^2$, set $\mathcal{H}_{-e}=(V\setminus e, \mathcal{E}\cap [V\setminus e]^2)$.
If for every $e,e'\in [V]^2$, $\mathcal{H}_{-e}\cong \mathcal{H}_{-e'}$, then we say that 
$\mathcal{H}$ is \emph{$-2$-monomorphic}. The following theorem shows that,  in the case of Steiner systems, join-irreducibility is equivalent to  $-2$-monomorphy.

\begin{theorem} \label{thm:steiner}Let $\mathcal{H}=(V, \mathcal{E})$ be a Steiner system. The following are equivalent:
 \begin{enumerate}[(i)]
  \item $f_{\mathcal{H}}$ is join-irreducible;
\item $\mathcal{H}_{ e}\cong \mathcal{H}_{ e'}$, for every $e,e'\in [V]^2$;
\item $\mathcal{H}$ is $-2$-monomorphic.
 \end{enumerate}
\end{theorem}

Note that for $k=2$ a Steiner system is a  complete graph. In this case $(i), (ii)$ and $(iii)$ hold simultaneously. For the proof of Theorem \ref{thm:steiner}, we suppose $k\geq 3$. We will need the following lemmas.

\begin{lemma}\label{lemma:0-Steiner}
Let $\mathcal{H}=(V, \mathcal{E})$ be a Steiner system. Then $D_{\mathcal{H}}= [V]^2$.
 \end{lemma}
\begin{proof} Observe first that  $\cup \mathcal{E}=V$. Next, let $e=\{i,j\}\in [V]^2$. Since $\mathcal{H}$ is a Steiner system and $k\geq 3$, a subset $E\subseteq V_e$ is an hyperedge of $\mathcal{H}_e$ if and only if $l_e\not \in E$ and $E$ is an hyperedge of $\mathcal{H}$, or $E=\{l_e\}\cup T\setminus e$, where $T$ is the unique hyperedge of $\mathcal{H}$ such that $e\subseteq T$, or $i\in T$ and $j\not \in T$ or $j\in T$ and $i\not \in T$. Since $\cup \mathcal{E}=V$, it follows that $\cup \mathcal{E}_e=V_e$, hence $e\in D_{\mathcal{H}}$. 
\end{proof}

\begin{lemma}\label{lemma:1-Steiner}
 Let $\mathcal H$ and $\mathcal H'$ be two Steiner systems and $e, e'\in [V]^2$. Then:
 \begin{enumerate}
 \item A map $f\colon V_{-e}\to V_{-e'}$ is an isomorphism from $\mathcal{H}_{-e}$ to $\mathcal{H}_{-e'}$ if and only if the map  $\overline{f}\colon V_{e}\to V_{e'}$ given by 
$\overline{f}(l_{e}) =l_{e'}$  and $\overline{f}(k)
=f(k)$ for $k\not =l_{e}$,  is an isomorphism from $\mathcal{H}_{e}$ to $\mathcal{H}_{e'}$.
\item \label{lemma:2-Steiner}
Every isomorphism $g\colon \mathcal{H}_{e} \to \mathcal{H}_{e'}$ is of the form $\overline{f}$ for some isomorphism 
$f\colon \mathcal{H}_{-e} \to \mathcal{H}_{-e'}$.
\end{enumerate}
\end{lemma}
\begin{proof} 
 (1) Clearly, $f$ is an isomorphism from $\mathcal{H}_{-e}$ to $\mathcal{H}_{-e'}$ whenever $\overline{f}$ is an isomorphism from $\mathcal{H}_{e}$ to $\mathcal{H}_{e'}$.

To show that the converse also holds, suppose that $f$ is an isomorphism.  Let
$T$ and $T'$ be the unique hyperedges of $\mathcal{H}$ containing $e$ and $e'$, respectively.

\begin{claim}\label{claim:2-Steiner}
 $f[T\setminus e]=T'\setminus e'$. In particular, 
$\overline{f}[\{l_e\}\cup (T\setminus e)]=\{l_{e'}\}\cup (T'\setminus e')$.
\end{claim}

\begin{proof}[Proof of claim] 
Let $i\in V\setminus e$;  denote by $d_{\mathcal{E}_{-e}}(i)$ the number of hyperedges in 
$\mathcal{E}_{-e}$ containing $i$ and define similarly $d_{\mathcal{E}}(i)$. Since $\mathcal{H}$ is a Steiner system, there are exactly one hyperedge containing  $i$ and intersecting $e$ if $i \in T$ (namely $T$) and  exactly two hyperedges when  $i\not \in T$. In other words,  the difference  $\eta(e, i):=d_{\mathcal{E}}(i)-d_{\mathcal{E}_{-e}}(i)$  is $1$ if $i\in T$ and $2$ if $i\in V\setminus T$. Furthermore $d_{\mathcal{E}}(i)=\frac{n-1}{k-1}$. Hence, with obvious notations $\eta(e, i)=\eta(e', f(i))$. From this,  $i\in T$ if  and only if $f(i) \in T'$ and thus the claim follows.
\end{proof}

Now, let  $S$ be an hyperedge of  $\mathcal{H}$ such that $S\cap e$ is a singleton.
\begin{claim}\label{claim:3-Steiner}
There is an hyperedge $S'$ of  $\mathcal{H}$ such that $S'\cap e'$ is a singleton and such that 
 $f[S\setminus e]=S'\setminus e'$. In particular, $\overline{f}[\{l_e\}\cup (S\setminus e)]=\{l_{e'}\}\cup (S'\setminus e')$.
\end{claim}

\begin{proof}[Proof of claim] 
Since $k\geq 3$, $\vert S\setminus e\vert \geq 2$. Let $j, j'\in S\setminus e$, $j\not =j'$ and let $S'$ be the unique hyperedge of 
$\mathcal{H}$ containing $\{f(j),f(j')\}$. Since $f$ is an isomorphism, we have $S'\cap e'\neq \emptyset$. By the previous claim we also have $S'\cap e'\neq e'$, and thus $S'\cap e'$ is a singleton.

Let $i\in S\setminus (e\cup \{j,j'\})$. We need to show that $f(i)\in S'$.
For that, let $S'_j$ and $S'_{j'}$ be the hyperedges of $\mathcal{H}$ containing 
$\{f(j),f(i)\}$ and $\{f(j'),f(i)\}$, respectively. Note first that $S', S'_j$ and $S'_{j'}$ intersect pairwise over $V\setminus e$. Moreover, since $f$ is an isomorphism, replacing $j,j'$ by $i,j$  and $S'$ by $S'_j$ we obtain, via   the  argument above,  that   $S'_{j}\cap e'$ is a singleton and, by the same token,  that  $S'_{j'}\cap e'$ is a  singleton too. Since this also holds for $S'\cap e'$, and $\vert e\vert =2$, two members of  $S', S'_{j}$ and $S_{j'}$ contain the same element of $e$.  Since $\mathcal H$ is a Steiner system, these two members coincide. Thus they contain $f(j)$ and $f(j')$.  Again from the fact that  $\mathcal H$ is a Steiner system, we have that   they coincide with $S'$. Hence $f(i)\in S'$ and the proof of the  claim is complete. 
\end{proof}
These two claim ensure that $\overline f$ is an isomorphism.

(2) Clearly, the lemma holds whenever $V$ is itself an hyperedge of $\mathcal{H}$. Thus, we may assume that this is not the case.
Now to prove the lemma, it is enough to show $g(l_e)=l_{e'}$ since in this case the restriction $f=g\mid_{V_{-e}}$ constitutes an isomorphism from $\mathcal{H}_{-e}$ to $\mathcal{H}_{-e'}$.

Let $T$  be the unique hyperedges of $\mathcal{H}$ containing $e$ and let   $\overline{T}={l_e}\cup (T\setminus e)$. Define  $T'$ similarly and let $\overline{T'}={l_{e'}}\cup (T'\setminus e')$.  Clearly  $\overline{T}$ is the only edge of $\mathcal H_{e}$ having size $k-1$. Since $g$ is an isomorphism of $\mathcal H_{e}$ on $\mathcal H_{e'}$, it maps  $\overline T$ on $\overline {T'}$. In particular,  $e'\in T'$.
Now observe that for each $i\not \in T$ there are exactly two hyperedges of $\mathcal H_{e}$ containing both $k$ and $l_e$, whereas for each $k\in T$ there is exactly  one, namely, $\overline{T}.$ Thus $g(l_e)=l_{e'}$.
\end{proof} 

\begin{proof}[Proof of Theorem \ref{thm:steiner}]

Implication  $(i)\Rightarrow(ii)$ follows from Lemma \ref {lemma:0-Steiner}. Implication $(ii)\Rightarrow (i)$ follows  from Proposition~\ref{cor:1}.
The implications $(ii)\Rightarrow (iii)$ and $(iii)\Rightarrow (ii)$ follow  respectively from (2) and (1) of Lemma~\ref{lemma:1-Steiner}. The proof of the theorem is now complete.
\end{proof}

\begin{problem}
For a Steiner triple system $\mathcal{H}=(V, \mathcal{E})$ does the following hold:
$\mathcal{H}$ is $-2$-monomorphic if and only if $Aut(\mathcal{H})$ is $2$-set transitive?
\end{problem}
Note that the automorphism group of a Steiner  systems is flag-transitive whenever it is $2$-set transitive. The converse holds for Steiner triple systems. There are several deep results about Steiner systems with a $2$-transitive or a flag transitive automorphism group (see the survey by Kantor \cite{kantor}). For example, any Steiner triple system with a $2$-transitive automorphism group must be a projective space over $GF(2)$ or an affine space over $GF(3)$, see e.g. \cite{hall, key}. The notion of monomorphy (with some of its variations) is due to R. Fra{\"\i}ss\'e. His book \cite{fraissetr} contains some important results concerning this notion. 

\section{Join-irreducible graphs}\label{irreduciblegraphs}

In this section, we  answer  Problem \ref{mainproblem} in the particular case of functions which are determined by undirected graphs, possibly with loops, that is, graphs $\mathcal{G} = (V, \mathcal{E})$ where $\mathcal{E}\subseteq [{V}]^2\cup [V]^1$.  As in the case of  hypergraphs, if we remove the  isolated vertices of  $\mathcal G$, the resulting graph $\check{\mathcal G}$ yields an equivalent function. In the sequel, when we speak of a \emph{join-irreducible} graph we simply mean a graph $\mathcal{G}$ such that $f_{\mathcal{G}}$ is join-irreducible. Note that  $\mathcal G$ is join-irreducible if and only if $\check {\mathcal G}$ is join-irreducible; note also that a join-irreducible graph must  satisfy $\card{\check{V}}\geq 2$.

Given a graph $\mathcal{G} = (V, \mathcal{E})$, we write $i\sim j$ if $\{i,j\}\in \mathcal{E}$. Set $V(i)=\{j\in V: i\sim j\}\cup \{i\}$. The \emph{degree} of a vertex $i$, denoted $d(i)$, is the cardinality $\card{V(i)}-1$. For example, in the \emph{complete} graph $K_n$ each vertex has degree $n-1$, while in a \emph{cycle} $C_n$ each vertex has degree $2$. Note that \emph{loops}, i.e., singleton edges, do not contribute to the degree of vertices. 

A graph $\mathcal{G}$ is said to be \emph{connected} if any two vertices of $\mathcal{G}$ are connected by a path. We denote by $\overline{\mathcal{G}}$ the \emph{complement}
of $\mathcal{G}$, that is, $\overline{\mathcal{G}} = (V, ([V]^2\setminus \mathcal{E})\cup ([V]^1\cap \mathcal{E}))$. 
The \emph{disjoint union} of two graphs $\mathcal{G}_1 = (V_1, \mathcal{E}_1)$, $\mathcal{G}_2 = (V_2, \mathcal{E}_2)$, $V_1\cap V_2=\emptyset $, is defined as the graph $\mathcal{G} = (V_1\cup V_2,\mathcal{E}_1 \cup \mathcal{E}_2)$. The \emph{graph join} of $\mathcal{G}_1 = (V_1, \mathcal{E}_1)$ and $\mathcal{G}_2 = (V_2, \mathcal{E}_2)$, denoted $\mathcal{G}_1+\mathcal{G}_2$, is defined as the disjoint union of $\mathcal{G}_1$ and $\mathcal{G}_2$ together with the edges $\{i_1,i_2\}$ for all $i_1\in V_1$ and $i_2\in V_2$. For further background in graph theory, see e.g.  \cite{bondy-murty, diestel, GGL}.

\subsection{Join-irreducible graphs: the loopless case}

In this subsection, we present an explicit description of those join-irreducible \emph{loopless} graphs, i.e., join-irreducible graphs
 $\mathcal{G} = (V, \mathcal{E})$ where $\mathcal{E}\subseteq [{V}]^2$. 
So throughout this subsection, we assume that $\mathcal{G}$ is a loopless graph. Furthermore, since $\mathcal G$ is join-irreducible if and only if $\check {\mathcal G}$ is join-irreducible, we also assume that $\mathcal G=\check {\mathcal G}$. We start with the disconnected case.

\begin{proposition}\label{disconnected}
 Suppose that $\mathcal{G}$ is disconnected. Then $\mathcal{G}$ is join-irreducible if and only if 
$\check{\mathcal{G}}$ is isomorphic to the disjoint union of $n$ copies of $K_3$, for some $n\geq 2$.
\end{proposition}
\begin{proof}
 Clearly, if ${\mathcal{G}}$ is isomorphic to the disjoint union of $n$ copies of $K_3$, for some $n\geq 2$, then $\mathcal{G}$ is join-irreducible. For the converse, let $\mathcal{G}_1$ and $\mathcal{G}_2$ be two connected components of $\mathcal{G}$. Note that $\card{\mathcal{G}_1},\card{\mathcal{G}_2}\geq 2$. Take $i\in \mathcal{G}_1$ and $j\in \mathcal{G}_2$, and let $e=\{i,j\}$. Clearly, $\card{\mathcal{G}_e}=\card{\mathcal{G}}-1$, no vertice is isolated in $\mathcal{G}_e$ and $\mathcal{G}_e$
has one less connected component than $\mathcal{G}$.

Now take $i,i'\in \mathcal{G}_1$ and let $e'=\{i,i'\}$. Clearly, for every such choice of $e'$, we have $e\not \approx e'$. Since $\mathcal{G}$ is join-irreducible, Theorem~\ref{theorem:key} implies that $\ess f_{\mathcal{G}_{e'}}<\ess f_{\mathcal{G}_e}$. In other words, for every $e'=\{i,i'\}$, $i,i'\in \mathcal{G}_1$, $\ess f_{\mathcal{G}_{e'}}\leq \ess f_{\mathcal{G}}-2$. From Theorem \ref{booleangap} it follows that $\mathcal{G}_1$ must be isomorphic to $K_3$. Since the choice of connected components was arbitrary, we conclude that $\mathcal{G}$ is isomorphic to the disjoint union of $n$ copies of $K_3$, for some $n\geq 2$. 
\end{proof}

To deal with the case of connected (loopless) graphs, we need to introduce some terminology.
Let $\mathcal{G}=(V, \mathcal{E})$ be a graph. A subset $S\subseteq V$ is said to be \emph{autonomous} if for every $i,i'\in S$ and $j\in V\setminus S$, $i\sim j$ if and only if $i'\sim j$. Moreover, $S$ is said
to be \emph{independent} if for every $i,i'\in S$, $i\not \sim i'$. For simplicity, we refer to autonomous independent sets as $ai$-sets.
We say that $\mathcal{G}$ is \emph{$ai$-prime} if its $ai$-sets are empty or singletons.

\begin{fact} For each $i\in V$, the union of all $ai$-sets containing $i$ is an $ai$-set called the \emph{$ai$-component} of $i$. Moreover, each graph $\mathcal{G}$ decomposes into $ai$-components.
\end{fact}

On the set of $ai$-components of $\mathcal{G}$ there is a graph structure, denoted $\mathcal{G}_{ai}$, in such a way that $\mathcal{G}$ is the lexicographic sum of its $ai$-components and indexed 
by $\mathcal{G}_{ai}$. Note that the graph $\mathcal{G}_{ai}$ is $ai$-prime. 

These constructions are variants of the classical notions of decomposition of graphs and prime graphs (see \cite{ehrenfeucht}). The following technical lemma is proved and extended to arbitrary graphs (not necessarily loopless) in the next subsection (see Lemma~\ref{lemma:auxGen}). 
\begin{lemma}\label{lemma:aux}
 Let $\mathcal{G}=(V, \mathcal{E})$ be a connected graph and suppose that there is $e \in [V]^2\setminus \mathcal{E}$ such that $\mathcal{G}_e$ has no isolated vertices. Then there is $e' \in \mathcal{E}$ such that $\mathcal{G}_{e'}$ has no isolated vertices.
\end{lemma}

Thus, if $\mathcal{G}=(V, \mathcal{E})$ is join-irreducible,  $\mathcal{G}_e$ has an isolated vertex for every $e\in [V]^2\setminus \mathcal{E}$. Moreover,  a nonedge  $e=\{i_1,i_2\}  \in [V]^2\setminus \mathcal{E}$ must be in a $ai$-component or there is $j\in V$ such that $d(j)=2$ and $i_1\sim j\sim i_2$.

We say that a graph $\mathcal{G}=(V, \mathcal{E})$ \emph{satisfies} $(P)$ if for every nonedge  $e=\{i_1,i_2\}  \in [V]^2\setminus \mathcal{E}$ there is $j\in V$ such that $d(j)=2$ and $i_1\sim j\sim i_2$.

Lemma \ref{lemma:aux} and the observation above yield the following.
\begin{corollary}  If  a connected graph $\mathcal G$ is join-irreducible, then $\mathcal G_{ai}$ satisfies $(P)$.
\end{corollary}

Our next proposition describes those graphs satisfying (P).

\begin{proposition}\label{prop:propP} A graph $\mathcal{G}=(V, \mathcal{E})$ satisfies $(P)$ if and only if $\mathcal{G}$ is either isomorphic to $K_n$, for some $n\geq 2$, $C_5$, $C_4$ or to a $3$-element path.
\end{proposition}

\begin{proof} We observe that each member of the list satisfies $(P)$. Conversely, suppose that $\mathcal{G}$ satisfies $(P)$. 
  
\begin{claim}\label{claim: P}Let $e=\{i_1,i_2\}  \in [V]^2\setminus \mathcal{E}$ and  $j\in V$ such that $d(j)=2$ and $i_1\sim j\sim i_2$. If $e':=\{j, j'\}\in [V]^2\setminus \mathcal{E}$ then either $e_1:=\{i_1,j'\}  \in \mathcal{E}$ and $d(i_1)=2$ or $e_2:=\{i_2,j'\}  \in \mathcal{E}$ and $d(i_2)=2$.
\end{claim}
\begin{proof}[Proof of Claim~\ref{claim: P}]
 Let $k \in V$ be such that $j\sim k\sim j'$. Since $V(j)=\{i_1,i_2\}$, we have $k=i_1$ or $k=i_2$. By property (P), it follows that either $e_1:=\{i_1,j'\}  \in \mathcal{E}$ and $d(i_1)=2$ or $e_2:=\{i_2,j'\}  \in \mathcal{E}$ and $d(i_2)=2$. 
\end{proof}

\begin{claim}\label{claim:clique} Let $i\in V$. If $d(i)\geq 3$, then $\mathcal{G}(i)=(V(i),[V(i)]^2\cap \mathcal{E})$ is isomorphic to $K_n$, for some $n\geq 3$. Moreover,
$\mathcal{G}(i)=\mathcal{G}$.
\end{claim}
\begin{proof}[Proof of Claim~\ref{claim:clique}] Let $i\in V$ be such that $d(i)\geq 3$ and, for the sake of contradiction, suppose that $\mathcal{G}(i)$ is not isomorphic to $K_n$, for some $n\geq 3$. Take $i_1,i_2\in V(i)$ such that $i_1\not \sim i_2$. By property (P), there is $j\in V$ such that $i_1\sim j\sim i_2$ and $d(j)=2$. Since $d(i)\geq 3$ and $d(j)=2$, $j\not \in V(i)$.
Let $j'\in V(i)\setminus \{i,i_1,i_2\}$ (such a $j'$ exists because $d(i)\geq 3$). Note that $j'\not\sim j$.
By Claim~\ref{claim: P}, either $e_1:=\{i_1,j'\}  \in \mathcal{E}$ and $d(i_1)=2$ or $e_2:=\{i_2,j'\}  \in \mathcal{E}$ and $d(i_2)=2$. This yields the desired contradiction because both $V(i_1)$ and $V(i_2)$ contain $\{i,j,j'\}$.  

To see that the last claim holds, suppose that there is $k\in V\setminus V(i)$.
Take $j\in V$ such that $i\sim j\sim k$. Since $\mathcal{G}(i)$ isomorphic to $K_n$, for some $n\geq 3$,
we have that $d(j)\geq 3$, which contradicts property (P). 
\end{proof}

According to  Claim \ref{claim:clique}, if $\mathcal G$ is not isomorphic to $K_n$, then the degree of each vertex is at most $2$. Since $\mathcal G$ satisfies $(P)$, it must be isomorphic to $C_5$, $C_4$ or to a $3$-element path.
\end{proof}

As a corollary we get the following result. 

\begin{corollary}\label{cor-ai}
If $\mathcal{G}$ is connected and 
join-irreducible, then $\mathcal{G}_{ai}$ is isomorphic to
 $K_n$, for some $n\geq 2$, or  to $C_5$.
\end{corollary}

Clearly,  each $K_n$, $n\geq 2$, and $C_5$ are join-irreducible graphs.
Thus, if $\mathcal{G}=(V, \mathcal{E})$ is an $ai$-prime graph, then $\mathcal{G}$ is join-irreducible if and only if ${\mathcal{G}}$ is isomorphic to $K_n$, for some $n\geq 2$, or to $C_5$.

Now if a connected and join-irreducible graph $\mathcal{G}=(V, \mathcal{E})$ is  not an $ai$-prime graph, then ${\mathcal{G}_{ai}}$ cannot be isomorphic to $C_5$. Indeed, for the sake of contradiction, suppose that ${\mathcal{G}_{ai}}$ is isomorphic to $C_5$. Let $\mathcal{G}_1, \ldots , \mathcal{G}_5$
be the $ai$-components of $\mathcal{G}$ such that $\mathcal{G}_i$ is connected to $\mathcal{G}_{i+1}$, for $i=1,2,3,4$ and $\mathcal{G}_5$ is connected to $\mathcal{G}_{1}$. Assume, without loss of generality, that $\card{\mathcal{G}_1}\geq 2$. Consider $i,i'\in \mathcal{G}_1$, $i_2\in \mathcal{G}_2$ and $i_3\in \mathcal{G}_3$, and let $e=\{i,i_2\}$ and $e'=\{i',i_3\}$. Clearly, $e\not \approx e'$ and $\ess f_{\mathcal{G}_{e}}=\ess f_{\mathcal{G}_{e'}}$. By Theorem~\ref{theorem:key} it follows that $\mathcal{G}$ is not  join-irreducible which constitutes the desired contradiction. 
Thus, by Corollary \ref{cor-ai} it follows that, in the non $ai$-prime case, if $\mathcal{G}=(V, \mathcal{E})$ is join-irreducible, then ${\mathcal{G}_{ai}}$ is isomorphic to $K_n$, for some $n\geq 2$.

\begin{proposition}Let $\mathcal{G}=(V, \mathcal{E})$ be a connected graph which is not $ai$-prime. Then $\mathcal{G}$ is join-irreducible if and only if ${\mathcal G}$ is isomorphic to one of the following graphs:
\begin{enumerate}[(i)]
\item ${K}_2+\overline{K}_m$,  for some $m\geq 2$;
\item $\overline{K}_n+\overline{K}_m$,
for some $n,m$ with $1\leq n<m$;
\item a graph join $\overline{K}_n+\ldots +\overline{K}_n$ of $r$ copies of $K_n$, for some $r, n\geq 2$.
\end{enumerate}
\end{proposition}
\begin{proof} As observed if $\mathcal{G}$ is join-irreducible, then ${\mathcal{G}_{ai}}$ is isomorphic to $K_r$, for some $r\geq 2$.  If $r=2$, then it is clear that  $\mathcal G$ is isomorphic to $\overline{K}_n+\overline{K}_m$ for some $n, m\geq 1$. 

Now suppose that $r=3$. If there is $i\in V$ such that $d(i)=2$, then $\mathcal G$ is isomorphic to ${K}_2+\overline{K}_m$,  for some $m\geq 2$. To verify the latter claim, let $\mathcal{G}_1,\mathcal{G}_2$ and $\mathcal{G}_3$ be the $ai$-components of $\mathcal G$ and suppose that there is $i\in V$ such that $d(i)=2$, say, $i\in C_1$.
Then, $\mathcal{G}_2$ and $\mathcal{G}_3$ are singletons, and since $\mathcal{G}$ is non $ai$-prime and ${\mathcal{G}_{ai}}$ is isomorphic to $K_3$, it follows that $\mathcal G$ is isomorphic to ${K}_2+\overline{K}_m$,  for some $m\geq 2$. If for every $i\in V$ we have $d(i)>2$, then $\mathcal G$ is isomorphic to $\overline{K}_{m_1}+\overline{K}_{m_2}+\overline{K}_{m_3}$, where at most one $ai$-component is a singleton.
Since $\mathcal{G}$ is join-irreducible, it follows that $m_1=m_2=m_3=n$, for some $n\geq 2$.

If $r > 3$, then for every $i\in V$ we have $d(i)>2$. Again from join-irreducibility, it follows that $\mathcal G$ is isomorphic to a graph join $\overline{K}_n+\ldots +\overline{K}_n$  of $r$ copies of ${K}_n$ for some $n\geq 2$.
\end{proof}

From these results, we obtain the description of the join-irreducible (loopless) graphs.

\begin{theorem} \label{thm:jigraphs}Let $\mathcal{G}=(V, \mathcal{E})$ be a (loopless) graph. Then $\mathcal{G}$ is join-irreducible if and only if $\check{\mathcal G}$ is isomorphic to one of the following graphs:
\begin{enumerate}[(i)]
 \item a disjoint union of $n$ copies of $K_3$, for some $n\geq 2$;
\item $C_5$;
\item ${K}_2+\overline{K}_m$,  for some $m\geq 2$;
\item $K_n$, for some $n\geq 2$;
\item $\overline{K}_n+\overline{K}_m$,
for some $n,m$ with $1\leq n<m$;
\item a graph join $\overline{K}_n+\ldots +\overline{K}_n$ of $r$ copies of $K_n$, for some $r, n\geq 2$.
\end{enumerate}
\end{theorem}

\subsection{Join-irreducible graphs: the general case}

In this subsection, we extend Theorem~\ref{thm:jigraphs} to arbitrary  (not necessarily loopless) graphs.  
The key result is the following. For a graph $\mathcal{G}=(V,\mathcal{E})$, let  $\mathcal{G}^0=(V,\mathcal{E}^0)$, where $\mathcal{E}^0=\mathcal{E}\cap [V]^2$. 

\begin{proposition}\label{prop:key}
Let $\mathcal{G}=(V,\mathcal{E})$ be any graph such that $V=\bigcup \mathcal{E}$. If $\mathcal{G}$ is join-irreducible, then $\mathcal{G}$ is either 
\begin{enumerate}
\item a disjoint union of loops or the disjoint union of loops and a copy of $K_3$, or
\item $\mathcal{G}^0$ is join-irreducible with no isolated vertices.
\end{enumerate}
\end{proposition}

To prove Proposition~\ref{prop:key}, we need to first extend Lemma~\ref{lemma:aux} to arbitrary graphs.

\begin{lemma}\label{lemma:auxGen}
Let $\mathcal{G}=(V,\mathcal{E})$ be a connected graph. Suppose that there are $i,j\in V$ such that $e=\{i,j\}\not \in \mathcal{E}$ and $\ess f_{\mathcal{G}_e}=\ess f_{\mathcal{G}}-1$. Then there are  $i',j'\in V$ such that $e'=\{i',j'\} \in \mathcal{E}$ and $\ess f_{\mathcal{G}_{e'}}=\ess f_{\mathcal{G}}-1$.
\end{lemma}

\begin{proof}
 Let $e=\{i,j\}\in [V]^2\setminus \mathcal{E}$ such that $\ess f_{\mathcal{G}_e}=\ess f_{\mathcal{G}}-1$. Since $\mathcal{G}$ is connected, there exists $i'\in V$ such that $i'\sim i$ or $i'\sim j$. Without loss of generality, assume that $d=\{i',i\}\in \mathcal{E}$.
 If $\ess f_{\mathcal{G}_{d}}=\ess f_{\mathcal{G}}-1$, then we are done. 
So suppose that $\ess f_{\mathcal{G}_{d}}<\ess f_{\mathcal{G}}-1$. 

If either $\{i\}, \{i'\}\in \mathcal{E}$ or $\{i\}, \{i'\}\not \in \mathcal{E}$, then there is $j'\in V$ such that $i'\sim j'\sim i$ and $d(j')=2$. Since $\mathcal{G}$ is connected, we have $d(i)>2$ or $d(i')>2$.
Suppose that $d(i)>2$. If $\{i\}, \{i'\}\in \mathcal{E}$, then $\ess f_{\mathcal{G}_{d'}}=\ess f_{\mathcal{G}}-1$, where $d'=\{i,j'\}$. 
If $\{i\}, \{i'\}\not \in \mathcal{E}$, then $\ess f_{\mathcal{G}_{e'}}=\ess f_{\mathcal{G}}-1$, where $e'=\{i',j'\}$. 
Similarly, the claim holds when $d(i')>2$. 

So we may assume that, say, $\{i\}\in \mathcal{E}$ and $\{i'\}\not \in \mathcal{E}$. The case $\{i'\}\in \mathcal{E}$ and $\{i\}\not \in \mathcal{E}$ can be verified similarly. 
If $l_d$, where $d=\{i,i'\}$, appears in some edge (singleton or pair) of $\mathcal{G}_d$, then there is $j'\in V$ such that $i'\sim j'\sim i$ and $d(j')=2$, and we have that $\ess f_{\mathcal{G}_{d'}}=\ess f_{\mathcal{G}}-1$, where $d'=\{i,j'\}$, or $\ess f_{\mathcal{G}_{e'}}=\ess f_{\mathcal{G}}-1$, where $e'=\{i',j'\}$, according to whether $d(i)>2$
or $d(i')>2$, respectively.

If $l_d$ does not appear in any edge of $\mathcal{G}_d$, then for every $k\in V$ we have that
$c=\{i, k\}\in \mathcal{E}$ if and only if $c'=\{i', k\}\in \mathcal{E}$. 
By connectivity, there is at least one such $k$. If there is exactly one, then $\ess f_{\mathcal{G}_{c'}}=\ess f_{\mathcal{G}}-1$ because $d(k)>2$. If there are at least two, then by choosing such a $k$ of greatest degree, it follows that 
 $f_{\mathcal{G}_{c}}=\ess f_{\mathcal{G}}-1$, and the proof of the lemma is complete. 
\end{proof}

In the sequel, we will also need the fact below.

\begin{fact}\label{fact:key}
 Let $\mathcal{G}=(V,\mathcal{E})$ be a connected loopless graph and let $i,i'\in V$.
Suppose that $\{l_e\}$, where $e=\{i,i'\}$, is not a member of $\mathcal{E}_e$. If $\ess f_{\mathcal{G}_{e}}=\ess f_{\mathcal{G}}-1$, then for any $\mathcal{G}^+=(V,\mathcal{E}^+)$, where $\mathcal{E}^+=\mathcal{E}\cup \{\{j\}\colon j \in S\}$ for some $S\subseteq V$, we have $\ess f_{{\mathcal{G}^+}_{e}}=\ess f_{\mathcal{G}^+}-1$.
\end{fact}

Finally, we will make use of the following result which shows that in most cases, reducibility can be verified using an edge and a nonedge.

\begin{lemma}\label{lemma:nonedge}
  Let $\mathcal{G}=(V,\mathcal{E})$ be a join-reducible, connected loopless graph. If for every $e=\{i,i'\}\in [V]^2\setminus \mathcal{E}$, $\mathcal{G}_{e}$ has at least one isolated vertex, then $\mathcal{G}$ is isomorphic to the graph join
$$\overline{K}_{m_1}+\overline{K}_{m_2}+\ldots +\overline{K}_{m_n},$$ 
where $n\geq 3$, $0<m_1\leq \ldots \leq m_n$ with $m_1<m_n$, and for $n=3$, $n_2\neq 1$.
\end{lemma}

\begin{proof}
  Necessarily, there is 
$e=\{i,i'\}\in [V]^2\setminus \mathcal{E}$, for otherwise $\mathcal{G}$ is a complete graph and thus it is join-irreducible.
Suppose first that there is such a pair $e=\{i,i'\}$ for which $l_e$ belongs to some member of $\mathcal{E}_e$. Then there is $k\in V$ such that $i\sim k\sim i'$ and $d(k)=2$. Since $\mathcal{G}$ is join-reducible, $\card{V}>3$ and hence there is $j\in V$ such that either $j\sim i$ but $j\not \sim i'$ or $j\sim i'$ but $j\not \sim i$. 

Suppose $j\sim i$ but $j\not \sim i'$ and let $d=\{k,j\}$. Since $\ess f_{\mathcal{G}_{d}}<\ess f_{\mathcal{G}}-1$, we must have $d(i)=2$. Also, there must exist $l$ such that $i'\sim l\sim j$ and $d(l)=2$, otherwise $\ess f_{\mathcal{G}_{c}}=\ess f_{\mathcal{G}}-1$, for $c=\{i',j\}$.
Moreover, $\ess f_{\mathcal{G}_{c}}<\ess f_{\mathcal{G}}-1$, for $c=\{i,l\},\{k,l\}$, and hence, $d(i')=d(j)=2$. Since $\mathcal{G}$ connected, it is isomorphic to $C_5$ which is a contradiction.
Thus, there is no $j\in V$ such that $j\sim i$ but $j\not \sim i'$.
Similarly, we can verify that there is no $j\in V$ such that $j\sim i'$ but $j\not \sim i$.

So we may ssume that, for every $e=\{i,i'\}\in [V]^2\setminus \mathcal{E}$, $l_e$ does not belong to any member of $\mathcal{E}_e$. Thus, for every $j\in V$ we have $j\sim i$ if and only if $j\sim i'$, i.e., $i$ and $i'$ belong to the same $ai$-component. Moreover, $\mathcal{G}_{ai}$ is isomorphic to a complete graph $K_n$. Since $\mathcal{G}$ is not join-irreducible, it has the form given in  the lemma.\end{proof}

\begin{proof}[Proof of Proposition~\ref{prop:key}]
Suppose that $\mathcal{G}$ is disconnected. Then each connected component $\mathcal{G}'$ having at least two vertices is join-irreducible. If $\mathcal{G}$ has an isolated loop, then $\gap f_{\mathcal{G}'}=2$. Using Theorem \ref{booleangap}, it can be verified that 
$\mathcal{G}'$ is a $K_3$. Moreover, it is not difficult  to see that there is no more than one such connected component.
If  there is no isolated loop, then by reasoning as in the proof of Proposition \ref{disconnected}, we can be shown that
 $\mathcal{G}$ is isomorphic to a disjoint union of $n\geq 2$ copies
of $K^+_3=(\{1,2,3\}, [\{1,2,3\}]^2\cup \{\{j\}\colon j \in S\})$, for some $S\subseteq \{1,2,3\}$.
This case is considered in Proposition \ref{prop:discjigraphs} below. 

 Suppose now that $\mathcal{G}$ is connected and that 
$\mathcal{G}^0$ is not join-irreducible. If  $\mathcal{G}^0$ isomorphic to the graph join 
\begin{equation}\label{exception}
 \overline{K}_{m_1}+\overline{K}_{m_2}+\ldots +\overline{K}_{m_n}
\end{equation}
as in Lemma \ref{lemma:nonedge}, then it is easy to see that any graph obtained from $\mathcal{G}^0$ by adding singletons $\{i\}$, $i\in V$, to the set of edges of $\mathcal{G}^0$, is join-reducible, and thus $\mathcal{G}$ is join-reducible.

So suppose that $\mathcal{G}^0$ is a join-reducible graph nonisomorphic to (\ref{exception}). By Lemma~\ref{lemma:nonedge}, there is $e=\{i,i'\}\in [V]^2\setminus \mathcal{E}^0$ such that $\mathcal{G}^0_{e}$ has no isolated vertices, that is, $\ess f_{\mathcal{G}^0_{e}}=\ess f_{\mathcal{G}^0}-1$. Since $\{l_e\}$ is not a member of $\mathcal{E}^0_e$, it follows from Fact~\ref{fact:key} that $\ess f_{{\mathcal{G}}_{e}}=\ess f_{\mathcal{G}}-1$. 
By Lemma~\ref{lemma:auxGen}, there is $e'=\{j,j'\}\in \mathcal{E}$ such that $\ess f_{\mathcal{G}_{e'}}=\ess f_{\mathcal{G}}-1$. Since $e\not\in \mathcal{E}$, we have $e\not \approx e'$ and thus $\mathcal{G}$ is join-reducible.

\end{proof}

By Proposition \ref{prop:key}, to completely describe the join-irreducible graphs (possibly with loops) we only need to focus on those graphs $\mathcal{G}$ without isolated loops such that $\mathcal{G}^0$ is join-irreducible. The description of the latter graphs is given in Theorem \ref{thm:jigraphs}. In the sequel, we consider a graph $\mathcal{G}=(V,\mathcal{E})$ without isolated loops.
 
The following proposition shows that, among those with no isolated loops, the disconnected join-irreducible graphs are exactly those which are loopless and join-irreducible.

\begin{proposition}\label{prop:discjigraphs}
If $\mathcal{G}^0$ is a disjoint union of $n$ copies of $K_3$, for some $n\geq 2$, then $\mathcal{G}$ is join-irreducible if and only if 
$\mathcal{G}=\mathcal{G}^0$.
\end{proposition}

\begin{proof}
By Theorem \ref{thm:jigraphs}, the condition $\mathcal{G}=\mathcal{G}^0$ is sufficient.
To show that it is necessary, observe first that no connected component $\mathcal{G}_i$ of $\mathcal{G}$ has loops in each vertex, for otherwise, by taking $i_1,i_2 \in \mathcal{G}_i$ and $j_1\in \mathcal{G}_j$, where $\mathcal{G}_j$ is a connected component of $\mathcal{G}$ different to $\mathcal{G}_i$, we have for $e=\{i_1,i_2\}$ and $e'=\{i_1,j\}$,
$$\ess f_{\mathcal{G}_{e}}=\ess f_{\mathcal{G}_{e'}}=\ess f_{\mathcal{G}}-1$$
but $e\not \approx e'$ which contradicts join-irreducibility by Theorem \ref{theorem:key}. 

We claim that no connected component of $\mathcal{G}$ has a loop. Indeed, suppose for the sake of contradiction that there is a connected component $\mathcal{G}_i$ of $\mathcal{G}$ with a loop. As we observed, we can find $i_1,i_2 \in \mathcal{G}_i$ such that $\{i_1\}$ is a loop, but not $\{i_2\}$. Now take $j_1\in \mathcal{G}_j$, where $\mathcal{G}_j$ is a connected component of $\mathcal{G}$ different to $\mathcal{G}_i$. 
We have for $e=\{i_1,j\}$ and $e'=\{i_2,j\}$,
$$\ess f_{\mathcal{G}_{e}}=\ess f_{\mathcal{G}_{e'}}=\ess f_{\mathcal{G}}-1$$
but it is easy to verify that $e\not \approx e'$.  This  contitutes the desired contradiction.
Thus $\mathcal{G}$ is a loopless graph and $\mathcal{G}=\mathcal{G}^0$.  
\end{proof}
 
The following propositions provide explicit descriptions of the remaining join-irreducible graphs, i.e., those graphs $\mathcal{G}$ for which $\mathcal{G}^0$ is isomorphic to one of the loopless graphs given in $(ii)-(vi)$ of Theorem \ref{thm:jigraphs}.
  
\begin{proposition}\label{prop:C5}
If $\mathcal{G}^0=C_5$, then $\mathcal{G}$ is join-irreducible if and only if $\mathcal{G}=\mathcal{G}^0$.
\end{proposition}

\begin{proof}
By Theorem \ref{thm:jigraphs}, the condition $\mathcal{G}=\mathcal{G}^0$ is sufficient.
 Conversely, suppose for the sake of contradiction that $\{i\}$ is a loop in  $\mathcal{G}$. 
 Take distinct $i_1,i_2\in V$ such that $i_1\sim i\sim i_2$. As before, we have for $e=\{i_1,i_2\}$ and $e'=\{i_1,i\}$
$$\ess f_{\mathcal{G}_{e}}=\ess f_{\mathcal{G}_{e'}}=\ess f_{\mathcal{G}}-1$$
but $e\not \approx e'$, contradicting the join-irreducibility of $\mathcal{G}$.
\end{proof}

\begin{proposition}\label{prop:Kn}
If $\mathcal{G}^0=K_n$, for $n\geq 2$, then $\mathcal{G}$ is join-irreducible if and only if it has $0,1,n-1$ or $n$ loops.
\end{proposition}

\begin{proof}
 If $\mathcal{G}$ has either $0$ or $n$ loops, then its automorphism group is 2-set transitive and, by Corollary \ref{2transitive}, 
 it is join-irreducible. If $\mathcal{G}$ has only one loop $\{i\}$, then for every distinct $i_1,i_2\in V\setminus \{i\}$, we have for $e=\{i,i_1\}$ and $e'=\{i_1,i_2\}$,   
$$\ess f_{\mathcal{G}_{e}}<\ess f_{\mathcal{G}_{e'}}=\ess f_{\mathcal{G}}-1$$
and, by Theorem \ref{theorem:key}, $\mathcal{G}$ is join-irreducible. Similarly, it can be verified that if $\mathcal{G}$ has $n-1$ loops, then it is join-irreducible.

Now suppose for the sake of contradiction that $\mathcal{G}$ has $n-k$ loops for $2\leq k\leq n-2$. Take distinct $\{i_1\},\{i_2\}\in \mathcal{E}$ and $\{j_1\},\{j_2\}\in [V]^1\setminus \mathcal{E}$. It is easy to see that, for $e=\{i_1,i_2\}$ and $e'=\{j_1,j_2\}$ 
$$\ess f_{\mathcal{G}_{e}}=\ess f_{\mathcal{G}_{e'}}=\ess f_{\mathcal{G}}-1$$
but $e\not \approx e'$ since $\mathcal{G}_{e}$ has 2 more loops than $\mathcal{G}_{e'}$. 
Thus $\mathcal{G}$ is not join-irreducible which constitutes the desired contradiction. 
\end{proof}

\begin{proposition}\label{prop:Kn+Km}
If $\mathcal{G}^0=\overline{K}_n+\overline{K}_m$,  for some $1\leq n,m$, then $\mathcal{G}$ is join-irreducible if and only if 
$\overline{K}_n$ and $\overline{K}_m$ have either no loops or loops in every vertex. 
\end{proposition}

\begin{proof}
 It is easy to verify that the conditions suffice to guarantee that $\mathcal{G}$ is join-irreducible. 
Let us prove that the converse also holds. The case when $m\geq n=1$ is straightforward. So let $m,n\geq 2$ and, for the sake of a contradiction,
suppose that there exist $i_1$ and $i_2$ in $\overline{K}_n$ or $\overline{K}_m$ such that $\{i_1\}$ is
a loop but $\{i_2\}$ is not a loop. 
Take $j$ in $\overline{K}_m$ or $\overline{K}_n$, according to whether $i_1$ and $i_2$ in $\overline{K}_n$ or $i_1$ and $i_2$ in $\overline{K}_n$, respectively. Clearly, we have for $e=\{i_1,j\}$ and $e'=\{i_2,j\}$,
$$\ess f_{\mathcal{G}_{e}}=\ess f_{\mathcal{G}_{e'}}=\ess f_{\mathcal{G}}-1$$
but $e\not \approx e'$. Thus $\mathcal{G}$ is not join-irreducible which yields the desired contradiction.
\end{proof}

\begin{proposition}\label{prop:K2+Km}
If $\mathcal{G}^0={K}_2+\overline{K}_m$,  for some $m\geq 2$, then $\mathcal{G}$ is join-irreducible if and only if either 
$K_2$ has no loops or two loops and $\overline{K}_m$ has no loops, or $\overline{K}_m$ has $m$ loops and $K_2$ has exactly one loop.
\end{proposition}

\begin{proof}
 It is easy to verify that the conditions suffice to guarantee that $\mathcal{G}$ is join-irreducible. 
To prove the converse, suppose first that  $K_2$ has exactly one loop, say, $i_1$ without and $i_2$ with a loop, and for the sake of a contradiction, suppose that $\overline{K}_m$ has at least one vertex $j$ without a loop. It is easy to see that we have, for $e=\{i_1,j\}$ and $e'=\{i_2,j\}$, 
$$\ess f_{\mathcal{G}_{e}}=\ess f_{\mathcal{G}_{e'}}=\ess f_{\mathcal{G}}-1$$
but $e\not \approx e'$ and thus $\mathcal{G}$ is not join-irreducible which contitutes the desired contradiction.  

Now suppose that $K_2$ has no loops or two loops and let $i_1$ and $i_2$ be its vertices. If every vertex of $\overline{K}_m$ has a loop
then, for every vertex $j$ of $\overline{K}_m$, we have for $e=\{i_1,i_2\}$ and $e'=\{i_1,j\}$, 
$$\ess f_{\mathcal{G}_{e}}=\ess f_{\mathcal{G}_{e'}}=\ess f_{\mathcal{G}}-1$$
but $e\not \approx e'$.
If there are vertices $j_1$ and $j_2$ of $\overline{K}_m$ such that $j_1$ has a loop but $j_2$ has no loop, then we have
for $e=\{j_1, j_2\}$ and $e'=\{i_1,j_2\}$, 
$$\ess f_{\mathcal{G}_{e}}=\ess f_{\mathcal{G}_{e'}}=\ess f_{\mathcal{G}}-1$$
but $e\not \approx e'$. In either case, we have that $\mathcal{G}$ is not join-irreducible and the proof of the proposition is complete.
\end{proof}

\begin{proposition}\label{prop:Kn+Kn}
If $\mathcal{G}^0=\overline{K}_n+\ldots +\overline{K}_n$ of $r$ copies of $K_n$, for some $r\geq 3$ and $ n\geq 2$, then $\mathcal{G}$ is join-irreducible if and only if $\mathcal{G}$ has no loops or loops in each vertex.
\end{proposition}

\begin{proof}
 It is easy to verify that the conditions suffice to guarantee that $\mathcal{G}$ is join-irreducible. 
To prove the converse, observe first that if there exist $i_1$ and $i_2$ in a $\overline{K}_n$ such that $\{i_1\}$ is a loop but not $\{i_2\}$, then by taking $j$ in another $\overline{K}_n$, we have $\{i_1,j\}\not \approx \{i_2,j\}$, and thus $\mathcal{G}$ is not join-irreducible. Hence, we may assume that each $\overline{K}_n$ has either no loops or loops in each vertex.

Now, if there is one $\overline{K}_n$ with loops in each vertex, and another with no loops, then $\mathcal{G}$ is not join-irreducible. Indeed, by taking $i_1$ in the former $\overline{K}_n$, $i_2$ in the latter $\overline{K}_n$, and another vertex $k$ such that $i_1\sim k\sim i_2$, we have 
for $e=\{i_1,k\}$ and $e'=\{i_2,k\}$, 
$$\ess f_{\mathcal{G}_{e}}=\ess f_{\mathcal{G}_{e'}}=\ess f_{\mathcal{G}}-1$$
but $e\not \approx e'$. This completes the proof of Proposition \ref{prop:Kn+Kn}.
\end{proof}

\end{document}